\documentclass[a4paper, 12pt, leqno]{amsart}
\usepackage{mathrsfs}
\usepackage{amsmath}
\usepackage{amscd}
\usepackage{amssymb}
\usepackage{amsthm}

\numberwithin{equation}{section}

\def\BN{\mathbb{N}}

\def\BZ{{\mathbb{Z}}}

\def\ll{\underset L\leq}
\def\rl{\underset R\leq}

\def\lrl{\underset {LR}\leq}
\def\llr{\lrl}
\def\el{\underset L\sim}
\def\er{\underset R\sim}
\def\elr{\underset {LR}\sim}

\def\cala{\mathcal{A}}

\begin{document}
\title[Lusztig's $a$-function for Coxeter groups]
{Lusztig's $a$-function for Coxeter groups with complete graphs}

\author{ Nanhua XI}
\address{HUA Loo-Keng Key Laboratory of Mathematics and Institute of Mathematics\\
Chinese Academy of Sciences\\
Beijing, 100190\\
China } \email{nanhua@math.ac.cn}
\thanks{}

\begin{abstract}
We show that  Lusztig's $a$-function of a Coxeter group is bounded
if the Coxeter group has a complete graph (i.e. any two vertices are
joined) and the cardinalities of finite parabolic subgroups of the
Coxeter group have a common upper bound.
\end{abstract}
\maketitle \setcounter{section}{-1}
\section{Introduction}
\def\tgt{\tilde G_2}
\label{sec:Intro} Lusztig's $a$-function for a Coxeter group is
defined in [L2] and is a very useful tool for studying cells in
Coxeter groups and related topics such as representations of Hecke
algebras. For an affine Weyl group, in [L2] Lusztig showed that the
$a$-function is bounded by the length of the longest element of the
corresponding Weyl group. It might be true that for any Coxeter
group of finite rank the $a$-function is bounded  by the length of
the longest element of certain finite parabolic subgroups of the
Coxeter group. In this paper we first show that this property
implies that the Coxeter group has a   lowest two-sided cell
(Theorem 1.5). We then show that Lusztig's $a$-function of a Coxeter
group has this property (Theorem 2.1) if the Coxeter group has a
complete graph (i.e. any  two different simple reflections of the
Coxeter group are not commutative) and the cardinalities of finite
parabolic subgroups of the Coxeter group have a common upper bound.
 For Coxeter groups of rank 3, Peipei Zhou  showed an analogue result by using the approach of the paper. These facts support
part (iv) of Question 1.13 in [X2].

\section{Preliminaries}

\noindent{\bf 1.1.} Let $(W,S)$ be a Coxeter group.  We  use $l$ for
the length function and $\le$ for the Bruhat order of $W$. The
neutral element of $W$ will be denoted by $e$.

Let $q$ be an indeterminate. The Hecke algebra $H$ of $(W,S)$ is a
free $\cala=\BZ[q^{\frac12},q^{-\frac12}]$-module with a basis
$T_w,\ w\in W$ and the multiplication relations are
$(T_s-q)(T_s+1)=0$ if $s$ is in $S$, $T_wT_u=T_{wu}$ if
$l(wu)=l(w)+l(u)$.

\def\tt{\tilde T}

For any $w\in W$ set $\tilde T_w=q^{-\frac{l(w)}2}T_w$. For any
$w,u\in W$, write $$\tilde T_w\tilde T_u=\sum_{v\in
W}f_{w,u,v}\tilde T_v,\qquad f_{w,u,v}\cala.$$ The following fact is
known and implicit in [L2, 8.3].

\medskip

 \noindent{(a)} For any $w,u,v\in W$, $f_{w,u,v}\in\cala $ is a polynomial in
$q^{\frac12}-q^{-\frac12}$ with non-negative coefficients  and
$f_{w,u,v}=f_{u,v^{-1},w^{-1}}=f_{v^{-1},w,u^{-1}}$. Its degree is
less than or equal to min$\{l(w),l(u), l(v)\}$.

Proof.  Note that $f_{x,y,e}=0$ if $xy\ne e$ and $f_{x,x^{-1},e}=1$
for any $x,y\in W$. Then it is easy to verify
$$f_{w,u,v}f_{v,v^{-1},e}=f_{w,w^{-1},e}f_{u,v^{-1},w^{-1}}.$$
So we have $f_{w,u,v}=f_{u,v^{-1},w^{-1}}=f_{v^{-1},w,u^{-1}}$. It
is clear that $f_{w,u,v}$ is a polynomial in
$q^{\frac12}-q^{-\frac12}$ with non-negative coefficients and
deg$f_{w,u,v}$ is less than or euqal to min$\{l(w),l(u)\}$. The
second assertion follows.

\medskip

For any $w,u,v$ in $W$, we shall regard $f_{w,u,v}$ as a polynomial
in $\xi=q^{\frac12}-q^{-\frac12}$. The following  fact is noted by
Lusztig [L3, 1.1 (c)].

\medskip

\noindent(b) For any $w,u,v$ in $W$ we have
$f_{w,u,v}=f_{u^{-1},w^{-1},v^{-1}}.$

\medskip

\noindent{\bf Lemma 1.2.} Let $(W,S)$ be a Coxeter group and $I$ is
a subset of $S$. The following conditions are equivalent.

(a) The subgroup $W_I$ of $W$ generated by $I$ is finite.

(b) There exists an element $w$ of $W$ such that $sw\le w$ for all
$s$ in $I$.

(c) There exists an element $w$ of $W$ such that $w\le ws$ for all
$s$ in $I$.

\medskip

Proof. Clear.

\medskip

We set $L(w)=\{s\in S\,|\, sw\le w\}$ and $R(w)=\{s\in S\,|\, ws\le
w\}$ for any $w\in W$.

\medskip

\noindent{\bf Lemma 1.3.} Let $w$ be in $W$ and $I$ is a subset of
$L(w)$ (resp. $R(w)$). Then $l(w_Iw)+l(w_I)=l(w)$ (resp.
$l(ww_I)+l(w_I)=l(w)$), here $w_I$ is the longest element of $W_I$.

Proof. Clear.

\medskip

\noindent{\bf 1.4.} For any $y,w\in W$, let $P_{y,w}$ be the
Kazhdan-Lusztig polynomial. Then all the elements
$C_w=q^{-\frac{l(w)}2}\sum_{y\le w}P_{y,w}T_y$, $w\in W$, form a
Kazhdan-Lusztig basis of $H$. It is known that
$P_{y,w}=\mu(y,w)q^{\frac12(l(w)-l(y)-1)}$ +lower degree terms if
$y<w$ and $P_{w,w}=1$.

For any $w,u$ in $W$, Write $$C_wC_u=\sum_{v\in W}h_{w,u,v}C_v,\
h_{w,u,v}\in\cala.$$ Following [L2], for any $v\in W$ we define
$$a(v)=\max\{i\in\BN\,|\, i=\text{deg}h_{w,u,v}, \ w,u\in W\},$$
here the degree is in terms of $q^{\frac12}$. Since $h_{w,u,v}$ is a
polynomial in $q^{\frac12}+q^{-\frac12}$, we have $a(v)\ge 0$.

We are interested in the bound of the function $a:W\to\BN$. Clearly,
$a$ is bounded if $W$ is finite. The following fact is known (see
[L3]) and easy to verify.

\medskip

\noindent(a) The $a$-function is bounded by a constant $c$ if and
only if deg$f_{w,u,v}\le c$ for any $w,u,v\in W$.

Lusztig showed that for an affine Weyl group the $a$-function is
bounded by the length of the longest element of the corresponding
Weyl group. This fact is important in studying cells in affine Weyl
groups. One consequence is that an affine Weyl group has a lowest
two-sided cell [S1].  We will show that the boundness of
$a$-function is also interesting in general.

\medskip

 Assume now that the $a$-function is bounded and its maximal value is
 $c$. Let $w,u,v$ be elements in $W$. We shall regard $h_{w,u,v}$ as a
polynomial in $\eta=q^{\frac12}+q^{-\frac12}$. Following Lusztig
[L2],  write
$h_{w,u,v}=\gamma_{w,u,v}\eta^{a(v)}+\delta_{w,u,v}\eta^{a(v)-1}+$lower
degree terms. Then $\gamma_{w,u,v}$ and $\delta_{w,u,v}$ are
integers. Let $\Omega$ be the subset of $W$ consisting of all
elements $w$ with $a(w)=c$.

Assume that $v\in\Omega$. For $w,u\in W$, we have
$f_{w,u,v}=\gamma_{w,u,v}\xi^{c}+$lower degree terms. Using 1.1(a)
and 1.1 (b) we get

\medskip

\noindent(b)
$\gamma_{w,u,v}=\gamma_{u,v^{-1},w^{-1}}=\gamma_{v^{-1},w,u^{-1}}=\gamma_{u^{-1},w^{-1},v^{-1}}$
for $w,u,v\in \Omega$.

\medskip

\noindent(c) Let $w,u\in W$ and $v\in\Omega$. If $\gamma_{w,u,v}\ne
0$, then $w,u$ are in $\Omega$ and
$\gamma_{w,u,v}=\gamma_{u,v^{-1},w^{-1}}=\gamma_{v^{-1},w,u^{-1}}=\gamma_{u^{-1},w^{-1},v^{-1}}$
is positive.

\medskip

\def\ll{\underset {L}{\leq}}
\def\rl{\underset {R}{\leq}}
\def\lrl{\underset {LR}{\leq}}
\def\llr{\lrl}
\def\el{\underset {L}{\sim}}
\def\er{\underset {R}{\sim}}
\def\elr{\underset {LR}{\sim}}
\def\ds{\displaystyle\sum}

Since $h_{w,u,v}\ne 0$ implies that $w\rl v,\ u\ll v$, by (c) we
obtain

\medskip

\noindent(d) Let $w,u\in W$ and $v\in\Omega$. If $\gamma_{w,u,v}\ne
0$, then $w\er v$, $w\el u^{-1}$ and $u\el v$. In particular,
$w,u,v$ are in the same two-sided cell.

\medskip

\noindent{\bf Theorem 1.5.} Let $(W,S)$ be a Coxeter group. Assume
that the $a$-function is bounded   by the length of the longest
element $w_0$ of a finite parabolic subgroup $P$ of $W$. Then the
two-sided cell of $W$ containing $w_0$ is the  lowest two-sided cell
of $W$. Moreover, the lowest two-sided cell contains all elements
$w$ in $W$ with $a(w)=l(w_0)$.

Proof.  We first show that $x\lrl w_0$ for any $x\in W$. (We refer
to [KL] for the definitions of the preorders $\lrl, \rl, \ll$ and
the equivalences $\el$, $\elr$ on $W$.)

Let $x\in W$ be such that $l(xw_0)=l(x)-l(w_0)$. We first show that
$x$ and $w_0$ are in the same left cell. Clearly $x\ll w_0$. Let
$y=xw_0$. Then
$$\tt_{x^{-1}}\tt_x=\tt_{w_0}(\sum_{z\in
W}f_{y^{-1},y,z}\tt_z)\tt_{w_0}.$$ Since $f_{y^{-1},y,e}=1$,
$f_{w_0,w_0,w_0}$ has degree $l(w_0)$ as a polynomial in
$\xi=q^{\frac12}-q^{-\frac12}$ and $f_{w,u,v}$ has non-negative
coefficients as a polynomial in $\xi$ for any $w,u,v\in W$, by
1.4(a) we conclude that $f_{x^{-1},x,w_0}$ has degree $l(w_0)$. Thus
$h_{x^{-1},x,w_0}$ has degree $l(w_0)$ as a Laurent polynomial in
$q^{\frac12}$. In particular, $h_{x^{-1},x,w_0}$ is non-zero, so
$w_0\ll x$. Hence $x$ and $w_0$ are in the same left cell.

Now assume that $x$ is an arbitrary element in $W$. Clearly there
exists $w\in P$ such that $l(xw)=l(x)+l(w)$ and
$l(xww_0)=l(xw)-l(w_0)$. We then have $ w_0\el xw\rl x$. Hence, the
two-sided cell containing $w_0$ is the lowest one among the
two-sided cells of $W$ (with respect to the partial order $\lrl$ on
the set of two-sided cells of $W$).

Now we show that the lowest two-sided cell contain all elements $w$
in $W$ with $a(w)=l(w_0)$.

Assume that $a(w)=l(w_0)$. Then there exists $x,y\in W$ such that
$\gamma_{x,y,w}\ne 0$ and $x,y,w$ are in the same two-sided cell. By
1.4 (c), $a(x)=a(y)=l(w_0).$ Choose $u\in P$ such that
$l(yu)=l(y)+l(u)=l(yuw_0)+l(w_0)$. It is easy to see that
$l(wu)=l(w)+l(u)=l(wuw_0)+l(w_0)$. Since
$\tt_x\tt_{yu}=(\tt_x\tt_y)\tt_u$, we have  $\gamma_{x,yu,wu}\ge
\gamma_{x,y,w}$. Thus $x,yu,wu$ are in the same two-sided cell. But
we have seen that $yu$ and $w_0$ are in the same two-sided cell.

The theorem is proved.

\medskip

\noindent{\bf Corollary 1.6} Let $(W,S)$ be a Coxeter group. Assume
that the $a$-function is bounded   by the length of the longest
element $w_0$ of a finite parabolic subgroup $P$ of $W$. Then
$$\{x\in W\,|\, l(xw_0)=l(x)-l(w_0)\}$$ is a left cell of $W$.

Proof. It follows from the proof of Theorem 1.5.

{\it Remark.} For affine Weyl groups, this result is due to Lusztig
[L2].

\medskip

\noindent{\bf 1.7.} Let $(W,S)$ be a Coxeter group. Assume that the
$a$-function is bounded   by the length of the longest element $w_0$
of a finite parabolic subgroup $P$ of $W$.  Denote the left cell
containing $w_0$ by $\Gamma$. Then $\Gamma=\{w\in W\,|\,
l(w)=l(ww_0)+l(w_0)\}.$ Let $J_{\Gamma\cap\Gamma^{-1}}$ be the free
$\BZ$-module with a basis $\{t_w\,|\, w\in \Gamma\cap\Gamma^{-1}\}$.
Define $t_wt_u=\sum_{v\in \Gamma\cap\Gamma^{-1}}\gamma_{w,u,v}t_v$.
Then $J_{\Gamma\cap\Gamma^{-1}}$ is an associative ring with unit
$1=t_{w_0}$.

Let $\Omega$ be the subset of $W$ consisting of all elements $w$
with $a(w)=l(w_0)$. We can define $J_{\Omega}$  and the
multiplication in $J_{\Omega}$ similarly. The multiplication is
associative. However, $J_{\Omega}$ has no unit in general, since
$\Omega$ contains infinite left cells in general, as shown in [B,
Be], see also Proposition 3.2.

\medskip

\noindent{\bf 1.8.  Remark.} Keep the assumption of Theorem 1.5.
Motivated by the work of Shi [S1, S2], we give some conjectures.

It is likely that the lowest two-sided cell is exactly the set of
elements $w$ in $W$ with $a(w)=l(w_0)$. Further, it is  likely that
the lowest two-sided cell coincides with the set of elements of $W$
of the form $xwy$ such that $l(xwy)=l(x)+l(w)+l(y)$, $l(w)=l(w_0)$
and $w$ is the longest element of a finite parabolic subgroup of
$W$.

Let $D'$ be the set consisting of all elements $x\in W$ such that

(1) $x=wy$ for some $w$ in a finite parabolic subgroup of $W$ with
length $l(w_0)$ and $y\in W$ and $l(x)=l(w)+l(y)$,

(2) for any $s$ in $L(w)$, there are no $z,z',u\in W$ such that
$sx=zuz'$, $l(sx)=l(z)+l(u)+l(z')$ and $u$ is in a finite parabolic
subgroup of $W$ with length $l(w_0)$.

 For any $x\in D'$, let $\Gamma_x$ be the subset of $W$
consisting of all elements  $zx$ satisfying $l(zx)=l(z)+l(x)$. It is
likely that $\Gamma_x$ is a left cell in the lowest-sided cell of
$W$ and the map $x\to \Gamma_x$ is a bijection between the set $D'$
and the set of left cells in the lowest two-sided cell. Also, the
set $D=\{y^{-1}wy\,|\, wy\in D'\}$ should be the set of
distinguished involutions in the lowest two-sided cell, here $wy$
satisfies the above  (1) and (2). When the Coxeter graph of $W$ is
connected we also conjecture that the set $D$ is finite if and only
if $W$ is finite or is an affine Weyl group or $st$ has infinite
order for any different simple reflections $s,t\in S$.

Assume that $wy$ satisfies (1) and (2). Let $zw\in W$ be such that
$l(zw)=l(z)+l(w)$. Then we should have
$C_{zw}C_{wy}=h_{w,w,w}C_{zwy}$. Also we should have
$\mu(z'wy,zwy)=\mu(z'w,zw)$ if $l(z'w)=l(z')+l(w)$. For affine Weyl
groups, these equalities are true, see [X1, SX].

If $(W,S)$ is crystallographic, then the function $a$ is constant on
a two-sided cell [L2]. Since $a(w_0)=l(w_0)$ (see [L2]), we see that
the lowest two-sided cell is exactly the set $\{w\in W\,|\,
a(w)=l(w_0)\}$.

For an affine Weyl group $W$, thanks to [S1, S2], we know that (a)
the lowest two-sided cell of $W$ coincides with the set of elements
of $W$ of the form $xwy$ such that $l(xwy)=l(x)+l(w)+l(y)$,
$l(w)=l(w_0)$ and $w$ is the longest element of a finite parabolic
of $W$; (b) $D$ is the set of distinguished involutions in the
lowest two-sided cell.

 In section 3 we will show that the above conjectures are true
 for certain Coxeter groups with complete graphs.

\section{Coxeter groups with complete graphs}

\def\st{\stackrel}
\def\sc{\scriptstyle}

Throughout this section $(W,S)$ is a Coxeter group and any two
simple reflections in $S$ are not commutative. In another words, the
Coxeter graph of $(W,S)$ is a complete graph. Another main result of
this article is the following.

\medskip

 \noindent{\bf Theorem 2.1.} Let $(W,S)$ be a Coxeter group. Assume that any two
different simple reflections are not commutative and the
cardinalities of finite parabolic subgroups of $W$ have a common
upper bound. Then Lusztig's $a$-function on $W$ is bounded by the
length of the longest element of certain finite parabolic subgroups
of $W$.

The remaining of this section is devoted to a proof of the theorem.

\medskip

\noindent{\bf Lemma 2.2.} Let $r,s,t$ be simple reflections such
that the orders of $rs,rt, st$ are greater than 2. Then there is no
element $w$ in $W$ such that $w=w_1r=w_2st$ and
$l(w)=l(w_1)+1=l(w_2)+2$.

Proof. We use induction on $l(w)$. When $l(w)=0,1,2,3$, the lemma is
clear. Now assume that the lemma is true for $u$ with length
$l(w)-1$. Since $r,t\in R(w)$, by Lemma 1.2, we know that the
subgroup $W_{rt}$ of $W$ generated by $r,t$ is finite. Let $w_{rt}$
be the longest element in $W_{rt}$. By lemma 1.3,
$w=w_3w_{rt}=w_4trt$ for some $w_3,w_4\in W$ and
$l(w)=l(w_3)+l(w_{rt})$, $l(w)=l(w_4)+3$. So we get $w_4tr=w_2s$.
Clearly we have $l(w_4rt)=l(w)-1$=$l(w_4)+2=l(w_2)+1$. By induction
hypothesis,  $w_2s$ does not exist, hence $w$ does not exist. The
lemma is proved.

\medskip

\noindent{\bf Lemma 2.3.} Keep the assumption of Theorem 2.1. Let
$x\in W$ and $t_1t_2\cdots t_m$ ($m\ge 2)$ be a reduced expression
of an element in $W$. Assume that $xt_1\le x$, $xt_2\cdots
t_{m-1}t_m\le xt_2\cdots t_{m-1}$, and  $l(xt_2\cdots
t_{m-1})=l(x)+m-2$. If for any reduced expression $s_1s_2\cdots s_m$
of $t_1t_2\cdots t_m$ with $xs_1\le x$ we have $l(xs_2\cdots
s_{m-1})$=$l(x)+m-2$, then $t_1t_2\cdots t_m$ is in a finite
parabolic subgroup of $W$ generated by two simple reflections.

Proof. If $m=2$, by Lemma 1.2, the result is clear. Now assume that
$m\ge 3$. Let $s=t_{m-1}$, $t=t_m$, and $y=xt_2\cdots t_{m-1}$. Then
$s,t\in R(y)$.  By Lemma 1.3, $y=y_1 s^a(ts)^b$ and
$l(y)=l(y_1)+a+2b$, here $a=0$ or 1 and $s^a(ts)^b$ is the longest
element of the subgroup $W_{st}$ of $W$ generated by $s,t$. Write
$t_1\cdots t_{m-1}=t_1\cdots t_is^d(ts)^c$, $t_i\ne t,s$, $d=0$ or
1, $c\ge 0$ and $d+2c+i=m-1$. We understand that $t_1\cdots t_i$ is
the neutral element $e$ of $W$ if $i=0$.

We need show that $i=0$. Since $t_m=t$ and $t_1t_2\cdots t_m$ is a
reduced expression, we must have $d+2c<a+2b$.

Assume  $a+2b=d+2c+1$. If $i\ge 1$, then $t_1t_2\cdots
t_{m-1}t_m=t_1\cdots t_is^a(ts)^b$ and $R(xt_2\cdots
t_i)\cap\{s,t\}$ contains exactly one element, denoted by $r$. (We
understand that $t_2\cdots t_i=e$ if $i=1$.) Then $t_1t_2\cdots t_m$
has a reduced expression of form $t_1t_2\cdots t_ir\cdots$. Since
$xt_1\le x$ and $i\le m-2$, by the assumptions of the lemma, we know
that $xt_2\cdots t_ir$ has length $l(x)+i$, this contradicts that
$xt_2\cdots t_ir\le xt_2\cdots t_i$. So $i=0$ in this case.

Assume $a+2b>d+2c+1$, then $b>c$ since $0\le a,d\le 1$. We have
$y_1s^a(ts)^b=(xt_1)t_1t_2\cdots t_is^d(ts)^c$. So
$y_1s^a(ts)^{b-c}s^d=(xt_1)t_1t_2\cdots t_i$ and
$l(s^a(ts)^{b-c}s^d)\ge 2$. Then $y_1s^a(ts)^{b-c}s^d=y_3st$ or
$y_3ts$ for some $y_3$ with $l(y_1s^a(ts)^{b-c}s^d)=l(y_3)+2$. If
$i\ge 1$,  we  have $y_1s^a(ts)^{b-c}s^dt_i\le y_1s^a(ts)^{b-c}s^d$.
Since $t_i\ne s,t$, by Lemma 2.2, this is impossible. The
contradiction leads  $i=0$. That is, all $t_1,t_2,...,t_m$ are in
$\{s,t\}$. The lemma is proved.

\medskip

{\it Remark.} The lemma is not true  in general. For instance, let
$(W,S)$ be of type $A_3$, $s_1,s_2,s_3$ are simple reflections such
that $s_1s_3=s_3s_1$. Consider $x=t_1t_2t_3t_4=s_2s_1s_3s_2$.

\medskip

\noindent{\bf Lemma 2.4.} Let $x\in W$ and $t_1\cdots t_m\cdots t_n$
($1<m<n$) be a reduced expression of an element in $W$. Assume that
(1) $l(xt_2\cdots t_{m-1}t_{m+1}\cdots t_{n-1})$ has length
$l(x)+n-3$, (2) $xt_1\le x$, (3) $xt_2\cdots t_{m-1}t_m\le
xt_2\cdots t_{m-1}$,  and (4) $xt_2\cdots t_{m-1}t_{m+1}\cdots
t_{n-1}t_n\le xt_2\cdots t_{m-1}t_{m+1}\cdots t_{n-1}$. Further,
assume that $t_1t_2\cdots t_m$ (resp. $t_m\cdots t_n$) is in a
parabolic subgroup $P$ (resp. $Q$) of $W$ with rank 2. Then $P=Q$ is
finite and $n=m+1$. In particular, $t_1\cdots t_n$ is in a finite
parabolic subgroup of $W$ generated by two simple reflections.

Proof. Let $t_m=s$ and $t_{m-1}=r$. Then $R(xt_2\cdots t_{m-1})$
contains $r,s$. Since the graph of $W$ is complete, any parabolic
subgroup of $W$ generated by more than two simple reflections is
infinite, by Lemma 1.2 we know that $R(xt_2\cdots t_{m-1})$ is
exactly $\{r,s\}$ and $P=<r,s>$ (the subgroup of $W$ generated by
$r,s$) is finite. Assume that $Q$ is generated by $s,t$. Clearly
$t_{m+1}\ne r,s$, so $t_{m+1}=t$. Let $xt_2\cdots t_{m-1}=yt_m$.
Then $l(yt_m)=l(y)+1$ and $R(y)$ does not contain $s$. We must have
$t\in R(y)$. Otherwise, $R(y)\cap\{s,t\}$ is empty and $xt_2\cdots
t_{m-1}t_{m+1}\cdots t_{n}=yt_mt_{m+1}\cdots t_n$ has length
$l(x)+n-2$. It contradicts the assumption $xt_2\cdots
t_{m-1}t_{m+1}\cdots t_{n-1}t_{n}\le xt_2\cdots t_{m-1}t_{m+1}\cdots
t_{n-1}.$ Therefore \break $xt_2\cdots t_{m-1}$=$yt_m=y_1ts=y_2srs$
has length $l(y_1)+2=l(y_2)+3$. So $y_1t=y_2sr$ has length
$l(y_1)+1)=l(y_3)+2$. By Lemma 2.2 we must have $t=r$ and then
$n=m+1$. The lemma is proved.

\medskip

\noindent{\bf  Lemma 2.5.} Let $x,w,y$ be elements in $W$. Assume
that $w$ is in a parabolic subgroup
 generated by two simple reflections $r,s\in S$, $l(w)\ge 3$ and $r,s$ are not in $R(x)\cup L(y)$. Then
 $l(xwy)=l(x)+l(w)+l(y).$

Proof: By Lemma 2.2, $R(xw)=R(w)$. Let $t_1\cdots t_n$ be a reduced
expression of $y$. Assume that $l(xwt_1\cdots
t_{m-1})=l(x)+l(w)+m-1$, $xwt_1\cdots t_{m-1}t_m$ $\le xwt_1\cdots
t_{m-1}$, and $m\le n$ is minimal for all reduced expressions of
$y$. Then $m\ge 2$. By Lemma 2.3, there exists $t_0\in R(w)$ such
that $t_0t_1\cdots t_{m-1}$ is in the finite parabolic subgroup of
$W$ generated by $t_0,t_1$. Since $l(w)\ge 3$ and $r,s$ are not in
$R(x)\cup L(y)$, by Lemma 2.2, $R(xwt_0)$ does not contain $t_1\in
L(y)$. Thus $t_0t_1\cdots t_{m-1}$ is the longest element of the
parabolic subgroup of $W$ generated by $t_0,t_1$ and $R(xwt_1\cdots
t_{m-1})=\{t_0,t_1\}$. So $t_0t_1\cdots t_{m-1}=t_1\cdots t_m$. Thus
$t_0\in L(y)$. This contradicts that $r,s\not\in R(x)\cup L(y)$. The
lemma is proved.

\medskip

\noindent{\bf Corollary 2.6.} Let $r,s$ be simple reflections and
$x,y,z\in W$ such that $x=yrs$ with $l(x)=l(y)+2$, $R(x)=\{s\}$,
$R(yr)=\{r\}$, $r,s\not\in L(z)$. Then $l(xz)=l(x)+l(z)$.

Proof. It follows from the proof of the above lemma.

\medskip

\noindent{\bf Lemma 2.7.} Assume that $w,u$ are  elements of a
finite parabolic subgroup $P$ of $W$ generated by two simple
reflections
 Then
 deg$f_{w,u,v}\le l(v)$ for $v\in P$   and
 $f_{w,u,v}=0$ if $v\not\in P$. (Recall that $f_{w,u,v}$ is a
 polynomial in $q^{\frac12}-q^{-\frac12}$.)

Proof. The first assertion follows from 1.1 9a) and the second
assertion is clear.

\medskip

\noindent{\bf Lemma 2.8.} Let $r,s,t$ be simple reflections and
$x,y,z\in W$. Assume that $x=yrs$, $R(yr)=\{r,t\}$, $R(x)=\{s\}$,
$R(y)=\{t\}$. If $r,s\not\in L(z)$, then deg $f_{x,z,w}\le 1$ for
all $w$ in $W$.

Proof. If $l(xz)=l(x)+l(z)$, nothing needs to prove. Assume that
$l(xz)<l(x)+l(z)$. Let $t_1t_2\cdots t_n$ be a reduced expression of
$z$. Then we can find a positive integer $m$  such that $xt_1\cdots
t_{m-1}t_m\le xt_1\cdots t_{m-1}$. By assumptions of the lemma,
clearly we have $m\ge 2$.  We choose the reduced expression of $z$
so that $m$ is minimal in all possibilities. According to Lemma 2.3,
$st_1\cdots t_{m-1}$ is in the parabolic subgroup of $W$ generated
by $s,t_1$.

We claim that  $t_1= t$. Otherwise, since $r,s\not\in L(z)$, Lemma
2.2 implies that the element  $st_1\cdots t_{m-1}=t_1\cdots t_m$ is
the longest element of the subgroup $<s,t_1>$ of $W$ generated by
$s,t_1$. This contradicts that $s\not\in L(z)$.

 Let $y_1\in W$ be such that $x=y_1trts$. Then $l(x)=l(y_1)+4$.
By Lemma 2.2 we know that $R(y_1tr)=\{r\}$. So $tst_1\cdots t_{m-1}$
is the longest element $w_{st}$ in $<s,t>$. Then (recall that
$\xi=q^{\frac12}-q^{-\frac12}$)
$$\tt_x\tt_z=\xi\tt_{y_1trw_{st}}\tt_{t_{m+1}\cdots t_k}+\tt_{y_1trw_{st}t_m}\tt_{t_{m+1}\cdots
t_k}.$$ We must have $s,t\not\in L(t_{m+1}\cdots t_n)$. Otherwise
$t_1\cdots t_{m+1}$ is the longest element of $<s,t>$ and $s\in
L(z)$, which contradicts our assumptions.

Since $l(w_{st})\ge 3$, by Lemma 2.5, we have
$$\tt_{y_1trw_{st}}\tt_{t_{m+1}\cdots t_k}=\tt_{y_1trw_{st}t_{m+1}\cdots
t_k}.$$

If $l(w_{st}t_m)\ge 3,$ using Lemma 2.5, we get

$$\tt_{y_1trw_{st}t_m}\tt_{t_{m+1}\cdots t_k}=\tt_{y_1trw_{st}t_mt_{m+1}\cdots
t_k}.$$ We are done in this case.

Assume now  $l(w_{st}t_m)=2$, then $w_{st}=tst$, $m=2$, $t_1=t$,
$t_2=s$. So $y_1trw_{st}t_m=y_1trst.$ If the longest element
$w_{rt}$ of $<r,t>$ is at least 4, then $w_{rt}t$ has length at
least 3. Since $s,t\not\in L(t_3\cdots t_n)$ and $r,s\not\in
L(t_1\cdots t_n)$ we see that $r,t\not\in L(stt_3\cdots t_n)$. Write
 $y_1trw_{st}t_m=y_2w_{rt}tst$, then $l(y_2w_{rt}tst)=l(y_2)+l(w_{rt}t)+2$.  By Lemma 2.5, we
know that
$$\tt_{y_1trw_{st}t_m}\tt_{t_{m+1}\cdots t_k}=\tt_{y_1trw_{st}t_mt_{m+1}\cdots
t_k}.$$ We are done in this case.

Assume now $w_{st}=sts$ and $w_{rt}=rtr$. Let
$u=y_1trw_{st}t_m=y_1trst.$ Assume that $s_1\cdots s_{n-2}$ is a
reduced expression of $t_3\cdots t_n$ and  $us_1\cdots s_{i-1}s_i\le
us_1\cdots s_{i-1}$ and $i$ is minimal in all possibilities. Note
that $R(u)=\{t\}$. By Lemma 2.3, $ts_1\cdots s_{i-1}$ is in a
parabolic subgroup of $W$ of rank 2. Since $s,t\not\in L(t_3\cdots
t_n)$, we  have $i\ge 2$. Also we have $s_1=r$ and $R(ut)=\{r,s\}$.
Otherwise, Lemma 2.2 implies that $ts_1\cdots s_{i-1}=s_1\cdots s_i$
is the longest element in $<t,s_1>$, so $t\in L(t_3\cdots t_n)$, a
contradiction. Now we have $R(u)=\{t\}$, $R(ut)=\{r,s\}$,
$R(uts)=\{r\}$ and $s,t\not\in L(t_3\cdots t_n)$. So can use
induction on $l(z)$ to see the lemma is true in this case.

The lemma is proved.

\def\tt{\tilde T}

\medskip

\noindent{\bf  2.9.} Now we prove Theorem 2.1. Let $x,y\in W$ and
consider
$$\tilde T_x\tilde T_y=\sum_{z\in W}f_{x,y,z}\tilde T_z.$$
We will prove that deg$f_{x,y,z}\le a_0$, here $a_0$ is the maximal
number among the lengths of the longest elements of all finite
parabolic subgroups of $W$.  Let $t_1t_2\cdots t_k$ be a reduced
expression of $y$. We may assume that $xt_1\le x$, otherwise we
replace $x$ by $xt_1$. We may further assume that $xs_1\le x$ for
any reduced expression $s_1\cdots s_m$ of $y$.

we use induction on $k$. For $k=0,1$, the result is clear. Now
assume that $k>1$.

If $xt_2\cdots t_k$ has length $l(x)+k-1$, then we have
$$\tilde T_x\tilde T_y=\xi \tilde T_{xt_2\cdots t_{k}}+\tt_{xt_1}\tt_{t_1y},$$
where $\xi=q^{\frac12}-q^{-\frac12}$. Using induction hypothesis we
see the theorem is true in this case.

Now assume that   $xt_2\cdots t_{m-1}t_m\le xt_2\cdots t_{m-1}$ for
some $2\le m\le k$. We may require that $m'\ge m>1$ if $s_1s_2\cdots
s_k$ is another reduced expression of $y$ and $xs_2\cdots
s_{m'-1}s_{m'}\le xs_2\cdots s_{m'-1}$. By Lemma 2.3, $t_1\cdots
t_m$ is in the parabolic subgroup $P$ generated by $t_1,t_2$.

Let $x_1$ (resp. $y_1$) be the  element in the coset $xP$ (resp.
$Py$) with minimal length.  Let $u,v\in P$ be such that $x=x_1w$ and
$y=uy_1$. Then we have
$$\tt_x\tt_y=\sum_{v\in P}f_{w,u,v}\tt_{x_1v}\tt_{y_1}.$$
By Lemma 2.7, deg$f_{w,u,v}\le l(v)$ and $v\in P$ if $f_{w,u,v}\ne
0$. If $l(v)\ge 3$, by Lemma 2.5, $l(x_1vy_1)=l(x_1v)+l(y_1)$. Hence
$\tt_{x_1w}\tt_{y_1}=\tt_{x_1wy_1}$. If $l(v)=2$, using Corollary
2.6 and Lemma 2.8 we see that deg$f_{x_1v,y_1,z}\le 1$ for any $z$.
If $l(v)=0$, by induction hypothesis, we see that the degrees of
$f_{x_1,y_1,z}$ are not greater than $a_0$ for any $z\in W$.

Now consider the case $l(v)=1$. In this case $v$ is a simple
reflection. We have $l(x_1v)=l(x_1)+1$  and $l(vy_1)=l(y_1)+1<l(y)$
since $m\ge 2$. Applying induction hypothesis to the equality
$$\tt_{x_1v}\tt_{vy_1}=\xi\tt_{x_1v}\tt_{y_1}+\tt_{x_1}\tt_{y_1},$$
we see that deg$f_{x_1v,y_1,z}\le a_0-1$ for any $z\in W$.

Theorem 2.1 is proved.

\medskip

\noindent{\bf Corollary 2.10.} Keep the assumption of Theorem 2.1.
Let $a_0$ be the maximal number among the lengths of the longest
elements of all finite parabolic subgroups of $W$. Then $a(w)=a_0$
 if and only if $w=xuy$ for some $x,y\in W$ and $u$ being
the longest element of a finite parabolic subgroup and
$l(w)=l(x)+l(y)+l(u)$, $l(u)=a_0$.

This is clear from the proof of Theorem 2.1.

\medskip

\section{Some consequences I - the lowest two-sided cell}

In this section $(W,S)$ is a Coxeter group with complete graph and
the cardinalities of finite parabolic subgroups of $W$ have a common
upper bound. We discuss the lowest two-sided cell of $W$. Let $a_0$
be the maximal value of the lengths of the longest elements of
finite parabolic subgroups of $W$ and let $\Lambda$ be the set
consisting of all the longest elements of finite parabolic subgroups
of the maximal cardinality (which is $2a_0$). Let $D'$, $D$ and
$\Gamma_x$ ($x\in D')$  be as in subsection 1.8.

\medskip

\noindent{\bf Proposition 3.1.} Keep the assumptions and  notations
above. We have

\noindent (a) The lowest two-sided cell of $W$ coincides with the
set $\{w\in W\,|\, a(w)=a_0\}$. So for any $x$ in the lowest
two-sided cell, there exists $y,z\in W$ and $u\in \Lambda$ such that
$x=zuy$ and $l(x)=l(z)+l(u)+l(y)$.

\noindent (b) The map $x\to \Gamma_x$ defines a bijection between
the set $D'$ and the set of left cells in the lowest two-sided cell
$c_0$.

\noindent (c)  The set $D=\{y^{-1}uy\,|\, w\in\Lambda,\  uy\in D',\
l(uy)=l(u)+l(y) \}$ is the set of distinguished involutions in the
lowest two-sided cell.

\noindent (d) Let $z,z'\in W$ and $uy\in D'$ be such that $u\in
\lambda$, $l(zuy)=l(z)+l(u)+l(y)$ and $l(z'uy)=l(z')+l(u)+l(y).$
Then $C_{zu}C_{uy}=h_{u,u,u}C_{zuy}$ and
$\mu(z'uy,zuy)=\mu(z'u,zu)$.

Proof. Let $$\Omega=\{zuy\in W\,|\, x,z\in W,\ u\in\Lambda,\text{
and } l(x)=l(z)+l(u)+l(y)\}.$$ We claim that $\Omega$ is the lowest
two-sided cell. Since $\Lambda\subset\Omega$, it suffices to prove
that $\Omega$ is a two-sided cell.

Let $x\in \Omega$. Then there exist $y,z\in W$ and $u\in \Lambda$
such that $x=zuy$ and $l(x)=l(z)+l(u)+l(y)$.  It is no harm to
assume that $uz$ is in $D'$. By computing $\tt_{zu}\tt_{uy}$ we see
easily that $\gamma_{zu,uy,w}\ne 0$ if and only if $w=x$ and
$\gamma_{zu,uy,x}=1.$ This implies that
$C_{zu}C_{uy}=h_{u,u,u}C_{zuy}$. The first part of (d) is proved.

Let $w\in W$ and $w\elr x$. Then there exist $w=w_1,w_2,...,w_n=x$
such that $\mu(w_i,w_{i+1})\ne 0$ or $\mu(w_{i+1},w_i)\ne 0$, and
$L(w_i)\not\subset L(w_{i+1})$ or $R(w_i)\not\subset R(w_{i+1})$ for
all $i=1,2,...,n-1$. We show that all $w_i$ are in $\Omega$. It is
no harm to assume that $n=2$ and  $L(w)\not\subset L(x)$. Let $s$ be
the simple reflection in $L(w)-L(x)$.  Then $C_{w}$ appears in
$C_sC_x$ with coefficient $\mu(w,x)$. Using the identity
$C_{zu}C_{uy}=h_{u,u,u}C_{zuy}$, we see that there exists $z_1\in W$
such that $l(z_1u)=l(z_1)+l(u)$, $C_{z_1u}$ appears in $C_sC_{zu}$
 and $\gamma_{z_1u,uy,w}\ne 0$. We must have $w=z_1uy$ and $\mu(z_1u,zu)=\mu(w,x)$ or $\mu(zu,z_1u)=\mu(x,w)$.
 So $w\in \Omega$. Part (a) is proved.

 Also we showed that
  $\Gamma_{uy}=\{zuy\,|\,z\in W,\ l(zuz)=l(y)+l(u)+l(y)\}$ is a left cell of $W$. Let $u_1y_1\in D'$, $l(u_1y_1)=l(u_1)+l(y_1)$,
  $u_1$ has length $a_0$ and
  is the longest element of a finite parabolic subgroup of $W$. If $u_1y_1\in\Gamma_{uz}$,
  then $u_1y_1=zuy$ for some $z\in W$ and $l(zuy)=l(z)+l(u)+l(y)$.
  By the definition of $D'$ we see that $z=e$. So $u_1y_1=uy$. Part (b) is proved.

 Comparing
the coefficients of $\tilde T_e$ in both sides of the  equality
$C_{zu}C_{uy}=h_{u,u,u}C_{zuy}$, we see that the
$l(w)-$2deg$P_{e,zuy}-a(u)\ge 0$. Moreover,
$l(w)-$2deg$P_{e,zuy}-a(u)= 0$ if and only if $z=y^{-1}$. In this
case, the coefficient of the term $q^{l(y)}$ is 1. Part (c) is
proved.

Now we prove the second part of (d).   Let $E_y,F_y\in H$ be such
that $C_uF_y=C_{uy}$ and $E_yC_u=C_{y^{-1}u}$. Then
$C_{zu}F_y=C_{zuy}$ and $E_yC_{uz^{-1}}=C_{(zuy)^{-1}}$. Thus
$h_{uz^{-1},z'u,w}=h_{y^{-1}uz^{-1}, z'uy,y^{-1}wy}$. Assume that
$z'uy< zuy$. Comparing the coefficients of $\tilde T_e$ in both
sides of the equality
$C_{(z'uy)^{-1}}C_{zuy}=h_{{z'uy}^{-1},zuy,w}C_{w}$, as in  [SX,
2.2], we see that the second part of (d) is true.

The proposition is proved.

\medskip

\noindent{\bf Proposition 3.2.} Let $(W,S)$ be a Coxeter group with
complete Coxeter graph and the cardinalities of finite parabolic
subgroups of $W$ have a common upper bound. Assume the cardinality
of $S$ is greater than 2 and the order of $st$ is finite for some
simple reflections $s,t$ in $S$. Then the number of left cells in
the lowest two-sided cell of $W$ is finite if and only if $W$ is an
affine Weyl group of type $\tilde A_2$.

Proof. The if part is clear (see [L2]). Now assume that $W$ is not
of type $\tilde A_2$. Let $s,t$ be simple reflections such that the
order $st$ is finite and maximal in all possibilities. Let $w$ be
the longest element of the
 subgroup $<s,t>$ of $W$ generated by $s,t$. Then $w$ is   in
the lowest two-sided cell  of $W$. If $w$ has length at least 4,
then $w(rst)^k$ is in $D'$ (see 1.8 for the definition of $D'$) for
any positive integer $k$, here $r$ is a simple reflection in
$S-\{s,t\}$. By Proposition 3.1 (b), we know that the number of left
cells in the lowest two-sided cell of $W$ is infinite.

If $w$ has length 3, then either $|S|\ge 4$ or one of  $rs$, $rt$
has infinite order for $r\in S-\{s,t\}$ since $W$ is not of type
$\tilde A_2$ and the length of $w$ is maximal among the longest
elements of finite parabolic subgroups of $W$. In first case, we can
find  two different simple reflections $r,v$ in $S-\{s,t\}$. Then
$w(rvst)^k$ is in $D'$ for any positive integer $k$.  In second
case, let $r\in S$ be different from $s,t$. It is no harm to assume
that $rs$ has infinite order. Then $w(rs)^k$ is in $D'$ for any
positive integer $k$. By Proposition 3.1 (b), in both cases the
number of left cells in the lowest two-sided cell of $W$ is
infinite.

The proposition is proved.

\medskip

\section{Some consequences II - other results}

In this section  $(W,S)$ is a Coxeter group such that any two simple
reflections are not commutative, except other specifications are
given. We shall give some other consequences of Theorem 2.1.

In [L1], Lusztig showed that the elements in $W$ with unique reduced
expressions forma a two-sided cell of $W$. If the order of $st$ is
 $\infty$ for  any two different simple reflections
$s,t$ of the Coxeter group $(W,S)$, then $W$ has only two two-sided
cells: $\{e\},\ W-\{e\}$, see [L3].

\medskip

\noindent{\bf Proposition 4.1.} Let $m\ge 3$ be a positive integer.
Assume that the order of $st$ is either $m$ or $\infty$ for  any two
different simple reflections $s,t$ of the Coxeter group $(W,S)$ and
the order of some $st$ is $m$. Then $W$ has only three  two-sided
cells.

 Proof. If $w\in W$ has
different reduced expressions, then there exist simple reflections
$s,t$ in $W$ and $x,y\in W$ such that $st$ has order $m$ and
$w=xuy$, $l(w)=l(x)+l(u)+l(y)$, where $u$ is the longest element in
the subgroup of $W$ generated by $s,t$. By Theorem 2.1, $m$ is the
maximal value of the $a$-function on $W$. According to Proposition
3.1, $w$ is in the lowest two-sided cell of $W$. Therefore, $W$ has
only the following three two-sided cells: $\{e\},$ $\{$elements in
$W$ with unique reduced expression$\}$, $\{$elements in $W$ having
different reduced expressions$\}$. The proposition is proved.

\medskip

\noindent{\bf 4.2.} Assume that any two  simple reflections in $W$
are not commutative. Let $O$ be the set of isomorphism classes  of
finite parabolic subgroups of $W$ with rank 2. It is likely that the
number of the two-sided cells of $W$ is $|O|+2$, here $|O|$ denotes
the cardinality of $O$. Proposition 4.1 supports this conjecture.
Below we will see that when $W$ is crystallographic, the conjecture
is also true.  We first establish some lemmas.

\medskip

\noindent{\bf Lemma 4.3.} Let $P$ and $Q$ be two different finite
parabolic subgroups of $W$ with rank 2. Denote  their longest
elements by $w$ and $u$ respectively. Assume $l(w)\le l(u)$. Let
$x,y\in Q$ be such that $l(wx)=l(w)+l(x)$, $l(wx)=l(wxu)-l(u)$,
$l(ywx)=l(y)+l(wx)$ and $l(y)=l(wx)-l(u)-1$. Then $\mu(u,ywx)=1$,

Proof. The existence of $x$ is clear. Since $l(w)\ge 3$, by Lemma
2.5,  $y$ exists. Using the  formulas (2.2.c) and (2.3.g) in [KL] we
can prove this lemma by a direct computation.

\medskip

\noindent{\bf Corollary 4.4.} Let $P$ and $Q$ be two finite
parabolic subgroups of $W$ with rank 2. Denote  their longest
elements by $w$ and $u$ respectively. Then $u\lrl w$ if $l(u)\ge
l(w)$. In particular, $w$ and $u$ are  in the same two-sided cell if
$l(w)=l(u)$.

Proof.  Let $y,x$ be as in Lemma 4.3. Since $l(w)\le (u)$, we have
$l(y)<l(x)$ and $L(ywx)$ is a proper subset of $L(u)$.  Lemma 4.3
then implies that $u\ll ywx$. Clearly, $ywx\lrl w$. So $u\lrl w$.
The lemma is proved.

\medskip

\noindent{\bf Proposition 4.5.} Let $(W,S)$ be a crystallographic
Coxeter group with complete Coxeter graph and $O$ be the set of
isomorphism classes  of finite parabolic subgroups of $W$ with rank
2. Then the number of the two-sided cells of $W$ is $|O|+2$.

Proof. By Theorem 2.1, the maximal value of the $a$-function on $W$
is at most 6. For $i=3,4,6$, let $W_i$ be the set of elements $x$ of
$W$ with the properties below: (1) $x=zuy$ for some $z,y\in W$, $u$
has length $i$ and is the longest element of a parabolic subgroup of
$W$, (2) $l(x)=l(z)+l(u)+l(y)$.

We have two obvious two-sided cells: $\{e\}$ and $\{$ elements in
$W$ with unique reduced expression$\}$. We claim that $W_6,
W_4-W_4\cap W_6$, $W_3-W_3\cap(W_4\cup W_6)$ are two-sided cells
whenever they are not empty.

First assume that $W_6$ is not empty. According to Proposition 3.1,
$W_6$ is the lowest two-sided cell of $W$. We claim that
$W_4-W_4\cap W_6$ is a two-sided cell of $W$ if it is not empty.
Clearly, $a(x)\ge 4$ for any $x\in W_4$. From the argument for
Theorem 2.1 we see easily that $a(x)\le 4$ if $x$ is not in $W_6$.
Since $W$ is crystallographic, by [L3, Corollary 1.9], for $x\in
W_4-W_4\cap W_6$ we have $x\elr u$, here $x=zuy$ is as in (1). Using
Corollary 4.4 we know that $W_4-W_4\cap W_6$ is in a two-sided cell.

Let $w$ be in $W_3$ but not in $W_4\cup W_6$. According to [L3,
Proposition 1.4], $\gamma_{w^{-1},w,d}=1$, here $d$ is the
distinguished involution in the left cell containing $w$. Using the
positivity  for Kazhdan-Lusztig polynomials of $W$ and for
$h_{w,w',w''}$ with $w,w',w''$ in $W$, we see that the argument for
Theorem 2.1 implies that $v\in W_6$ if deg$h_{w^{-1},w,v}\ge 4$.
Since $a(w)=a(d)\ge 3$ and $d$ is not in $W_6$, we must have
$a(w)=3$ and $w$ is not in the two-sided cell containing $x$.
Therefore, $W_4-W_4\cap W_6$ is a two-sided cell if it is not empty.
We also showed that $W_3-W_3\cap(W_4\cup W_6)$ is a two-sided cell
if it is not empty.

If $W_6$ is empty, the discussion is similar and simpler. The
proposition is proved.

\section{Some comments}

In this section we propose two questions. Let $(W,S)$ be an
arbitrary Coxeter group.  In [L2], it is showed that $a(w)\le (w)$
for any $w$ in a Weyl group. This result was extended to arbitrary
crystallographic Coxeter groups by Springer, see [L3] for a proof.
It is natural to suggest that $a(w)\le l(w)$ for $w$ in an arbitrary
Coxeter group.

Assume that $(W,S)$ is connected (i.e. its Coxeter graph is
connected). Let $P$ and $Q$ be two finite parabolic subgroups of
$W$. It is likely that the longest elements of $P$ and $Q$ are in
the same two-sided cell of $W$ if $P$ and $Q$ are isomorphic Coxeter
groups.

\bibliographystyle{unsrt}

\begin{thebibliography}{}

\bibitem[B]{B} B\'edard, R.: {\sl Left V-cells for hyperbolic Coxeter groups,} Comm.
Alg. \textbf{17} (1989), no. 12, 2971-2997.

\bibitem[Be]{Be} Belolipetsky, M.: {\sl Cells and representations of right-angled
Coxeter groups,} Selecta Math. (N.S.) \textbf{10} (2004), no. 3,
325-339.

\bibitem[KL]{KL}
Kazhdan, D., Lusztig, G.: {\sl Representations of Coxeter groups and
Hecke algebras,}  Invent. Math. \textbf{53}  (1979), 165-184.

\bibitem[L1]{L1}Lusztig, G.: {\sl Some examples of square integrable representations of semisimple $p$-adic
groups,}
 Trans. Amer. Math. Soc. \textbf{277} (1983), no. 2, 623-653.

\bibitem[L2]{L:cell2} Lusztig, G.:  {\sl Cells in affine Weyl groups,} in ``Algebraic groups and
related topics", Advanced Studies in Pure Math., vol. \textbf{6},
Kinokuniya and North Holland, 1985, pp. 255-287.

\bibitem[L3]{L:cell3} Lusztig, G.:  {\sl Cells in affine Weyl groups,
II,} J. Alg. \textbf{109} (1987), 536-548.

\bibitem[SX]{SX}  Scott, L. and Xi, N.: {\sl Some non-trivial Kazhdan-Lusztig coefficients of an affine Weyl
group of type $\tilde A_n$,} Science China Mathematics \textbf{53}
(2010), No. 8, 1919-1930.

\bibitem[S1]{S1} Shi, J.: {\sl A two-sided cell in an affine Weyl group,}
 J. London Math. Soc. (2) \textbf{36} (1987), no. 3, 407-420.

\bibitem[S2]{S2} Shi, J.: {\sl A two-sided cell in an affine Weyl group,II,}
 J. London Math. Soc. (2) \textbf{37} (1988), no. 2, 253-264.

\bibitem[X1]{Xi1}  Xi, N.: The based ring of the lowest two-sided cell of
an affine Weyl group, J. Alg. \textbf{134} (1990), 356-368.

\bibitem[X2]{Xi2}
Xi, N.: Representations of affine Hecke algebras, Lecture Notes in
Mathematics, \textbf{1587}. Springer-Verlag, Berlin, 1994.

\bibitem[Z]{Z}
Zhou, Peipei: Lusztig's $a$-function for Coxeter groups of rank 3.
In preparation.
\end{thebibliography}

\end{document}